\documentclass[11pt]{article}

\title{Ivan Bernoulli Series Universalissima}
\author{A.K.Kwa\'sniewski\\  
\\ Institute of Computer Science, \textsl{Bia\l ystok University}(*)\\
PL - 15-887 Bialystok , ul.Sosnowa 64,  Poland
\\e-mail: kwandr@wp.pl
\\(*) former \textsl{Warsaw University Division}}

\usepackage{amsmath,amsthm}

\chardef\bslash=`\\ 
\begin{document}
\maketitle

\begin{abstract}
\noindent One states, recalls here that the textbooks formula
named as celebrated Taylor formula pertains historically to Johann
Bernoulli (1667 - 1748). When Ivan Bernoulli`s \textit{Series
Universalisima }was published in Acta Eruditorum in Leipzig in
1694 Brook Taylor was nine years old. As on the 01.01.2006 one
affirms 258-th anniversary of death of Johann Bernoulli this
\textbf{calendarium-article} is \textit{Pro}-posed \textit{Pro hac
vice}, \textit{Pro memoria}, \textit{ Pro nunc }and \textit{Pro
opportunitate} - as circumstances allow.
\end{abstract}

\vspace{2mm}


\vspace{1mm}

\noindent AMS S.C. (2000) 01A45, 01A50, 05A40, 81S99.

\vspace{2mm}

\section{ The first historical remarks}

\noindent \textbf{Principal Data}

\noindent\textbf{ 1.I.} Johann \textbf{Bernoulli} [1667 - 1748]-
\emph{"Archimedes of his age"}. Johann Bernoulli attained great
fame already in his lifetime. He was elected a fellow of the
academies of  St Petersburg, Paris, Berlin, London and Bologna. He
was called the \emph{"Archimedes of his age"} and this is what is
inscribed on his tombstone. His also famous as a Master of the
brilliant pupil Leonhard.

\noindent\textbf{ 1.II.} Leonhard \textbf{Euler} [1707-1783].
Leonhard Euler - also born in Switzerland - studied under Giovanni
Bernoulli at Basel. Euler and Bernoullies became also
Personalities of the Great Imperial St. Petersburg of those times.
Euler supported a lifelong friendship with Ivan Bernoulli's sons
Daniel. Jean or Giovanni or Johann is Ivan there in Russia`s
Capital St. Peterburg as in Sonin`s article title [1].

\noindent Note:  Johann Bernoulli`s  "Series universalissima" was
already at work in his Acta Eruditorum paper published in 1694 in
Leipzig.

\noindent Complete list is available via address below.

\noindent \textbf{Bibliography-Data:} \noindent Johann (Giovanni,
Jean, Ivan) Bernoulli (27.01.1667 - 01.01.1748)

\noindent List of publications of Johann Bernoulli including post
mortem papers with quotations from  and on Johann Bernoulli are
attainable via Brandenburgische Akademie der Wissenschaften
Akademiebibliothek  in
$$http://bibliothek.bbaw.de/kataloge/literaturnachweise/bernou-1/literatur.pdf$$

\noindent \textbf{1.III.} Brook\textbf{ Taylor }[1685 - 1731].

\noindent Note: Taylor's Methodus incrementorum directa et inversa
(1715) had given birth to what is now called the "calculus of
finite differences". It also contained the celebrated formula
known as Taylor's expansion, the importance of which was
recognized only in 1772 when Lagrange proclaimed it the basic
principle of the
differential calculus.\\

\noindent \textbf{ 1.IV.}  N.Y.\textbf{Sonin}.Nikolay Yakovlevich
Sonin [1849 - 1915] was the pupil of Nicolai Vasilievich Bugaev
and went on to make a major contribution to mathematics. He taught
at the \textbf{University of Warsaw} where he obtained a doctorate
in 1874. Then in 1876 N.Y. Sonin  was appointed to a chair in the
University of Warsaw. In 1894 Sonin moved to St Petersburg.
Together with A A Markov, Sonin prepared a two volume edition of
Chebyshev's works in French and Russian. Here for us of primary
importance is Sonin`s article on Ivan Bernoulli. The source of
detailed information and historian like investigation paper. This
is Academician at that time  Sonin`s article  [1] dated from
\textbf{1897}. There (p. 341) he quotes and supports the
conviction of an anonymous author of \emph{Epistola pro eminente
mathematico, Dn. Johanno Bernulio, contra quedam ex Anglia
antagonistam scripta} in Acta Eruditorum 1716 p. 307 - the
conviction that one deals here with plagiarism
in Taylor`s Methodus (1715).\\
According to  Salvatore Anastasio (see more in the third section)
"Historians have since shown that Taylor was not guilty of
plagiarism, but only of having failed  to keep up on the
literature from the Continent."

\noindent In the same spirit in one of their internet articles
$$http://www-groups.dcs.st-and.ac.uk/~history/Mathematicians/Taylor.html$$
 by J. J. O'Connor and E. F. Robertson  the authors indicate that
what we now call Taylor series  is the one which \emph{Taylor was
the first to discover}. James Gregory [ 1638 - 1675],  Newton
[1643-1724], Leibniz [1646-1716], Johann Bernoulli   [1667 - 1748]
and de Moivre [1667 - 1754] had all known \emph{variants of
Taylor's Theorem}. James Gregory, for example, knew that $arctan x
= x - x3/3 + x5/5 - x7/7 + ...$.  J. J. O'Connor and E. F.
Robertson state that "\emph{all of these mathematicians had made
their discoveries independently, and Taylor's work was also
independent of that of the others. The importance of Taylor's
Theorem remained unrecognised until 1772 when Lagrange proclaimed
it the basic principle of the differential calculus. The term
"Taylor's series" seems to have used for the first time by
Lhuilier in 1786}".

\noindent Well. Those were the times  without e-mails though
famous margins and covers were prosperously used...and here  again
the same authors  about James Gregory:\\
\emph{"In February 1671 he discovered Taylor's theorem (not
published by Taylor until 1715), and the theorem is contained in a
letter sent to Collins on 15 February 1671. The notes Gregory made
in discovering this result still exist written on the back of a
letter sent to Gregory on 30 January 1671 by an Edinburgh
bookseller. Collins wrote back to say that Newton had found a
similar result and Gregory decided to wait until Newton had
published before he went into print. He still felt badly about his
dispute with Huygens and he certainly did not wish to become
embroiled in a similar dispute with Newton"}.\\
\noindent  This quotation above comes  from
$$http://www-groups.dcs.st-and.ac.uk/~history/Mathematicians/Gregory.html$$
\noindent The similar story is with Graves and GHW algebra. Though
this is another story of the HiStory. History of Mathematics.

\vspace{1mm}

\noindent \textbf{1.V. C.Graves} ad (*) from Key Words see  [2] by
C.~Graves from (1853--1857). The Graves-Heisenberg-Weyl (GHW)
algebra is called by physicists the Heisenberg - Weyl algebra.
See: [3,4] and references therein). We shall deal with GHW algebra
in the sequel. Before that let us continue the main story.\\

\noindent \textbf{The first historical remark with nowadays form
of formulas.} \noindent Here are the famous examples of expansion
$$\partial_0=\sum_{n=1}^{\infty}\frac{x^{n-1}}{n!}\frac{d^n}{dx^n}$$
or  $$\epsilon_0=\sum_{n=0}^{\infty}(-1)^n
\frac{x^{n}}{n!}\frac{d^n}{dx^n}$$ where $\partial_0$ is the
divided difference operator while $\epsilon_0$ is at the zero
point evaluation functional. If one compares these with
\textit{"series universalissima"} of J.Bernoulli from \textit{Acta
Erudicorum} (1694) (see commentaries in \cite{3,4,5}) and with
$$exp\{yD\}=\sum_{k=0}^{\infty}\frac{y^kD^k}{k!},\ \
D=\frac{d}{dx},$$ then confrontation with B.Taylor's {\em
"Methodus incrementorum directa et inversa"} (1715), London;
entitles one to call the expansion formulas considered above and
in this note {\em "Bernoulli - Taylor formulas"} or (for $n
\rightarrow
\infty$) {\em "Bernoulli - Taylor series" }\cite{3,4,5}.\\
\noindent  In support to this conviction postulate - facts come to
become the support. \noindent \textbf{Facts.} In  1694   Johann
Bernoulli publishes articles with his series univeralissima
applied to several functions.

\noindent [J.Bernoulli 1] Modus generalis construendi omnes
Aequationes differentiales primigradus. Lipsiae 1694 in: Acta
Eruditorum. 1694. S.435-437

\noindent [J.Bernoulli 2] Additamentum effectionis omnium
Quadraturarum et Rectificationum Curvarum per seriem quandam
generalissimam. Lipsiae 1694 in: Acta Eruditorum. 1694. S.437-441

\noindent From Sonin [1] p. 341  we learn  more. In  August of
1695 in the letter to Leibniz  Johann Bernoulli among others -
considered  also the function $f(x) = x^x$ and he also calculated
its integral over $[0,1]$  interval using his series
universalissima. This was highly appraised by Leibniz who  quoted
this result without proof in his article "Principia Calculi
Exponentialium, seu Percurrentium" on mart of 1697 in Acta
eruditorum, Leipzig.

\noindent Well, after  then series universalissima was not unknown
and more . It was really  "universalissima"  in the operational
spirit of Leibniz approach to Calculi - the up today  modern and
efficient attitude to Calculi. It is so for example in
combinatorics with Laplace`s generating functions` rigorous
apparatus that had been harnessed to work efficiently since its
firmly established birth in Bernoulli`s Acta Eruditorum articles.

\vspace{2mm}

\section{ The second historical remark on Johann Bernoulli and Brook Taylor as both the Great}

\noindent  \textbf{Information 2.1.} Johann \textbf{Bernoulli} was
elected a fellow of the academy of St Petersburg.  Johann
Bernoulli - the Discoverer of \textit{Series Universalissima} was
"\textbf{Archimedes of his age}" and this is indeed inscribed on
his tombstone.

\noindent  \textbf{Information 2.2.}Brook \textbf{Taylor.}
"\textit{There are other important ideas which are contained in
the Methodus incrementorum directa et inversa of 1715 which were
not recognised as important at the time. These include singular
solutions to differential equations, a change of variables
formula, and a way of relating the derivative of a function to the
derivative of the inverse function. Also contained is a discussion
on vibrating strings, an interest which almost certainly come from
Taylor's early love of music.} "

\noindent  Brook "\textit{\textbf{Taylor }was a mathematician of
far greater depth than many have given him credit for:- A study of
Brook Taylor's life and work reveals that his contribution to the
development of mathematics was substantially greater than the
attachment of his name to one theorem would suggest. His work was
concise and hard to follow. The surprising number of major
concepts that he touched upon, initially developed, but failed to
elaborate further leads one to regret that health, family concerns
and sadness, or other unassessable factors, including wealth and
parental dominance, restricted the mathematically productive
portion of his relatively short life.}"

\noindent These extracts come from the article by J. J. O'Connor
and E. F. Robertson
$$http://www-groups.dcs.st-and.ac.uk/~history/Mathematicians/Taylor.html$$

\section{ The third historical remark with champagne is a quotation}

\noindent This third historical remark with champagne in Champagne
is a quotation of an unsurpassed feature short note by Salvatore
Anastasio on Bernoulli-Taylor litigation.
$$http://www.newpaltz.edu/math/fall98.pdf The End of a Feud   by$$
\noindent \textbf{Quotation} \noindent \textit{ "The End of a
Feud" } by Sal Anastasio \textbf{(1998)}.

\noindent Many students of mathematics are aware of the
controversy between Isaac Newton and Gottfried  Wilhem Leibniz
over which of them was the first to "discover" the Calculus. For
many years in the early eighteenth century, followers of Newton
and Leibniz disputed bitterly over the priority issue. However,
historians today are convinced that each man deserves credit as an
independent "discoverer". \noindent What is less well known is the
275 year long feud between the families of Brook Taylor (1685 -
1731) and Johann Bernoulli (1667 - 1748), only recently ended.
\noindent  Brook Taylor (yes, he of "Taylor Series" fame) was an
English mathematician deeply committed to the cause of Newton.
\textbf{Bernoulli}, equally committed to Leibniz, was a Swiss
mathematician who \textbf{greatly developed Leibniz' techniques
(which proved far superior to Newton's and are those common today)
} and was responsible for the first textbook in Calculus:
L'Hospital's Analyse des Inniment Petits, 1696.   ("L'Hospital's"
? Yes, but that's another story!) \noindent It appears that Taylor
and Bernoulli began their own dispute in 1715, when Taylor
published  results in England which had already been discovered on
the Continent by Leibniz and Bernoulli. \noindent (Historians have
since shown that Taylor was not guilty of plagiarism, but only of
having failed  to keep up on the literature from the Continent.)
\noindent The French probabilist, Pierre Remond de Montmort,
apparently the Jimmy Carter of his day, tried very hard, but in
vain, to patch things up between Taylor and Bernoulli. \noindent
But, not to worry. About 8 years ago, on July 7, 1990, Francois de
Montmort, a descendant of Pierre, played host, at the de Montmort
ancestral chateau in Champagne, to descendants of Bernoulliand
Taylor: Rene Bernoulli of Basel, Switzerland, and Chalmers Trench
of Slane, Ireland.

\noindent After toasting each other with champagne (in Champagne),
the two journeyed to the front lawn of  the chateau with a shovel
in tow. There they solemnly dug a hole (would I lie?) and
successfully buried a hatchet (I kid you not) which had been
provided by an American historian of science. The feud is now officially over."\\
\noindent \textbf{End of quotation.}

\noindent \section{ Closing and rather technical remarks - for
today} Remarks to be rather skip over by non specialists.

\noindent \textbf{Introductory information.} An extension of Ivan
Bernoulli formula  of a new sort with the rest term of the Cauchy
type was recently derived also by the author in the case of the so
called $\psi$-difference calculus which constitutes the
representative case of extended umbral calculus. Naturally the
central importance of such a type formulas is beyond any doubt -
and recent publications do confirm this historically established experience.\\

\noindent \textbf{4.I.} All formulas based on Bernoulli-Taylor
formula for extended umbral calculus  from [3] (and corresponding
references by the author therein) may be quite easily extended
\cite {5} to the case of any $Q\in End(P)$ linear operator that
reduces by one the degree of each polynomial \cite {6}. Namely one
introduces:
\textbf{Definition 4.1.}
    $$\hat {x}_{Q} \in End(P), \hat x_{Q}: F[x] \to F[x] $$
     such that  $(x^{n}) = \frac{{\left( {n
+ 1} \right)}}{{\left( {n + 1} \right)_{\psi} } }q_{n + 1}; n \geq
0;$ where $Qq_n=nq_{n-1}$.

\noindent Then $\star_Q$ product of formal series  and
$Q$-integration are defined almost mnemonic analogously [3,4,5].

\noindent \textbf{4.II.} In  1937   Jean  Delsarte  [7]  had
derived the general Bernoulli-Taylor formula for a  class of
linear operators $\delta$ including  linear operators that reduce
by one the degree of each polynomial. The rest term of the
Cauchy-like type in his Taylor formula (I) is given in terms of
the unique solution of a first order partial differential equation
in two real variables. This first order partial differential
equation  is determined by the choice of  the linear operator
$\delta$  and the function  f under expansion. In  our
Bernoulli-Taylor $\psi$-formula  or in its straightforward
$\star_Q$ product of formal series and $Q$-integration
generalization -  there is no need to solve any partial
differential equation.

\vspace{1mm}

\noindent  \textbf{4.III.}   In  [8] (1941) Professor  J. F.
Steffensen - the Master of polynomials application to actuarial
problems
$$(see:http://www.math.ku.dk/arkivet/jfsteff/stfarkiv.htm)$$
supplied a remarkable derivation of another Bernoulli-Taylor
formula with the rest of "Q-Cauchy type" in the example presenting
the "Abel poweroids"

\vspace{1mm}

\noindent  \textbf{4.IV.} The recent paper [9] (2003) by Mourad E.
H. Ismail, Denis Stanton may serve as a kind of indication  for
pursuing further investigation. There the authors have established
two new q-analogues of a Taylor series expansion for polynomials
using special Askey-Wilson polynomial bases. As "byproducts" their
important paper  includes also  new summation theorems,quadratic
transformations for q-series and new results on a q-exponential
function.

\vspace{1mm}

\noindent  \textbf{4.V.} Let us also draw an attention to two more
different publication on the subject which are the ones referred
to as [10, 11]. The $q$-Bernoulli theorems are named here and
there above as $q$-Taylor theorems. The corresponding
$(q,h)$-Bernoulli theorem for the $\partial _{q,h}$-difference
calculus of Hahn (see [12]) might be obtained the same way as the
$q$-Bernoulli-Taylor theorem constituting here the special case of
the Viskov method application. This is so because the $\partial
_{q,h}$-difference calculus of Hahn may be reduced to $q$-calculus
of {\bf Thomae-Jackson} [13] due to the following observation. Let
$$ h \in R ,(E_{q,h}\varphi)(x)=\varphi (qx+h)$$
and  let
\begin{equation}\label{}
(\partial _{q,h}\varphi)(x)= \frac{\varphi(x)-
\varphi(qx+h)}{(1-q)x-h}
\end{equation}
Then (see Hann in [12])
\begin{equation}\label{eq}
\partial_{q,h} = E_{1,{\frac{{-h}}{{1-q}}}} \partial_q E_{1,{\frac
{h}{{1-q}}}}.
\end{equation}
It is easy now to derive corresponding formulas including
Bernoulli-Taylor $\partial _{q,h}$-formula obtained  in [14] by
the Viskov method [12,3]  which for
$$q\rightarrow 1 , h\rightarrow 0 $$
recovers the content of one of Viskov examples while for
$$q \rightarrow 1 , h \rightarrow 1$$
one recovers the content of the another Viskov example. The case
$h \rightarrow 0 $ is included in the formulas of $q$-calculus of
Thomae-Jackson easy to be specified : see  [12] (see also
\textbf{thousands} of up-date references there). For Bernoulli-
Taylor Formula (presented during $PTM$ - Convention Lodz - $2002$)
: contact [13] for its recent version.

\vspace{1mm}

\noindent \textbf{4.VI.} As in [5]  the rule $\tilde{g}(x)=g(\hat
{x}_{\psi}){\bf 1}$ defines the map $\sim  : g \mapsto \tilde{g}$
which is an umbral operator   $ \sim : P\mapsto P $. It is
mnemonic extension of the corresponding $q$ - definition by
Kirschenhofer [14] and  Cigler [15]. This umbral operator (without
reference to [14, 15]) had been already used in theoretical
physics aiming at Quantum Mechanics on the lattice [16]. The
similar aim is represented by  [17,18] (see further references
there)  where incidence algebras are being prepared for that
purpose. As it is well known - the classical umbral [19, 20] and
extended [5] finite operator calculi  may be formulated in the
reduced incidence algebra language. Hence both applications of
related tools to the same goal are expected to meet at the arena
of GHW algebra description of both [14,13].

\vspace{1mm}

\noindent\textbf{ At the end} let us again come back to $17$-th
and $18$-th centuries springs of modern Calculi.

\vspace{1mm}

\noindent Recall. \textbf{Johann Bernoulli} (27.01.1667 -
01.01.1748). We thus affirm  at this very 2006 year the
\textbf{258-th} anniversary of death of Johann Bernoulli
\textbf{just today on 01.01.2006.}

\vspace{1mm}

\noindent \textbf{His Legacy Bibliography-Data:}

\noindent List of publications of Johann Bernoulli including post
mortem papers with quotations from  and on Johann Bernoulli is
attainable via Berlin-Brandenburgische Akademie der Wissenschaften
Akademiebibliothek  in
$$http://bibliothek.bbaw.de/kataloge/literaturnachweise/bernou-1/literatur.pdf$$

\end{document}